\documentclass[a4paper,11pt]{article}

\usepackage{graphicx}
\usepackage{amsfonts}
\usepackage{amsmath, amssymb}
\usepackage{amsthm}
\usepackage{psfrag}

\newtheorem{theorem}{Theorem}
\newtheorem{proposition}{Proposition}
\newtheorem{definition}{Definition}

\DeclareMathOperator{\Tr}{tr} \DeclareMathOperator{\Det}{det}

\begin{document}

\title{Families of traveling impulses and fronts in some models with cross-diffusion}

\author{Faina S Berezovskaya$^{1}$ , Artem S Novozhilov$^{2}$, Georgy P Karev$^{2,}$\footnote{Corresponding author: tel.: +1 (301) 451-6722; fax: +1 (301) 435 7793; e-mail: karev@ncbi.nlm.nih.gov}\\
\textit{\normalsize $^{1}$Howard University, 6-th Str. Washington DC 20059}\\
\textit{\normalsize $^{2}$National Institutes of Health, 8600
Rockville Pike Bethesda MD 20894} }

\date{}

\maketitle

\begin{abstract}

An analysis of traveling wave solutions of partial differential
equation (PDE) systems with cross-diffusion is presented. The
systems under study fall in a general class of the classical
Keller-Segel models to describe chemotaxis. The analysis is
conducted using the theory of the phase plane analysis of the
corresponding wave systems without a priory restrictions on the
boundary conditions of the initial PDE. Special attention is paid
to families of traveling wave solutions. Conditions for existence
of front-impulse, impulse-front, and front-front travelling wave
solutions are formulated. In particular, the simplest mathematical
model is presented that has an impulse-impulse solution; we also
show that a non-isolated singular point in the ordinary
differential equation (ODE) wave system implies existence of
free-boundary fronts. The results can be used for construction and
analysis of different mathematical models describing systems with
chemotaxis.

\end{abstract}

\section{Introduction}\label{intro}
In this paper we study one-dimensional traveling wave solutions
for the models in the form
\begin{equation}\label{basic}
    \begin{split}
U_t&=(\mu U_x-f(U,V)V_x)_x,\\
V_t&=g(U,V).
\end{split}
\end{equation}
Here $\mu>0$ is a constant; $f(U,V)$ and $g(U,V)$ are functions
whose properties will be specified later; $U=U(x,t)$, $V=V(x,t)$;
in the following we put $\mu=1$ without loss of generality.

The model \eqref{basic} is known as a particular case of the
classical Keller-Segel models to describe chemotaxis, the movement
of a population $U$ to a chemical signal $V$ (see, e.g.,
\cite{Horstmann,Keller1,Keller2}). In system \eqref{basic} $\mu$
denotes the constant diffusion coefficient; $f(U,V)/U$ is the
chemotactic sensitivity, which can be either positive or negative;
$g(U,V)$ describes production and degradation of the chemical
signal; it is customary to include also in the second equation of
\eqref{basic} the diffusion term of the form $D V_{xx}$ which
would describe diffusion of the chemical signal, but, we adopt
hereafter, that, to a first approximation, $D$ can be taken zero.
For biological interpretation of the solutions of \eqref{basic} we
refer to the cited literature, references therein, and to Section
3.2; macroscopic derivation of equation \eqref{basic} can be found
in, e.g., \cite{Erban}.

The chemotactic models are the partial differential equations
(PDEs) with cross-diffusion terms; these systems possess special
mathematical peculiarities \cite{Ni}. Such systems were used,
e.g., to model the movement of traveling bands of \textit{E. coli}
\cite{Keller1,Keller2,Bud}, amoeba clustering \cite{Mikhailov},
insect invasion in a forest \cite{Berez2}, species migration
\cite{Feltman}, tumor encapsulation and tissue invasion
\cite{Habib,Sherratt}. Many different spatially non-homogeneous
patterns can be observed in chemotactic models, for a survey see,
e.g., \cite{tsyganov} and references therein. One such pattern is
that of traveling waves which spread through the population.

A \textit{traveling wave} is a bounded solution of system
\eqref{basic} having the form
$$
U(x,t)=U(x+ct)\equiv u(z),\quad V(x,t)=V(x+ct)\equiv v(z),
$$
where $z=x+c\,t$ and $c$ is the speed of wave propagation along
$x$-axis; $u(z)$ and $v(z)$ are the wave profiles ($u$-profile and
$v$-profile respectively) of solution ($U(x,t),\,V(x,t)$).

Substituting these traveling wave forms into \eqref{basic} we
obtain
\begin{equation*}
    \begin{split}
cu'&=(u'-f(u,v)v')',\\
cv'&=g(u,v),
\end{split}
\end{equation*}
where primes denote differentiation with respect to $z$.

On integrating the first equation in the last system we have
\textit{the wave system} of \eqref{basic}:
\begin{equation}\label{ws}
    \begin{split}
    u' & =cu+f(u,s)g(u,s)/c+\alpha,\\
      v'  & =g(u,s)/c.
\end{split}
\end{equation}

Here $\alpha$ is the constant of integration that depends on the
boundary conditions for $U(x,t)$ and $V(x,t)$. In various
applications it is usually possible to determine this constant
prior to analysis of the wave system. For instance, considering
system \eqref{basic} as a model of chemotactic movement
\cite{Keller1}, where the variable $U(x,t)$ plays the role of the
population density and $V(x,t)$ is an attractant, one usually
supposes that $\int_{-\infty}^\infty u(x,t)\, dx$ should be
finite, which implies that $\alpha=0$ (e.g., \cite{Horstmann}). On
the contrary, in our analysis we do not specify the boundary
conditions for \eqref{basic} and consider $\alpha$ as a new
parameter.

Each traveling wave solution of \eqref{basic} has its counterpart
as a bounded orbit of \eqref{ws} for some $\alpha$; in our study
we elucidate the following question: for which $\alpha$ there
exist traveling wave solutions of \eqref{basic} and describe all
such solutions. We also note that the case of $\alpha=0$ does not
exclude a model with infinite mass of $U(x,t)$ if the traveling
wave solution is a front; moreover, the solutions corresponding to
finite mass can be only impulses (see below for the terminology).

It is worth noting that due to specific form of system
\eqref{basic} with cross-diffusion terms the wave system has the
same dimension as the initial system \eqref{basic}, which
significantly simplifies the analysis. This is one of
peculiarities which distinguishes cross-diffusion PDEs from those
with only diffusion terms (see also \cite{Berez2, Berez4}).

We shall study possible wave profiles of \eqref{basic} and their
bifurcations with changes of the parameters $c$ and $\alpha$ by
the methods of phase plane analysis and bifurcation theory
\cite{An,Kuznet}. In this way, the problem of describing all
traveling wave solutions of system \eqref{basic} is reduced to the
analysis of phase curves and bifurcations of solutions of the wave
system \eqref{ws} without a priory restrictions on boundary
conditions for \eqref{basic}.

There exists a known correspondence between the bounded traveling
wave solutions of the spatial model \eqref{basic} and the orbits
$u(z),\,v(z)$ of the wave system \eqref{ws} (e.g., \cite{Berez2,
Berez, Murray,Volpert}) that we only list for the cases most
important for our exposition.

\begin{figure}[thb]
  \centering
  \includegraphics[width=0.8\textwidth]{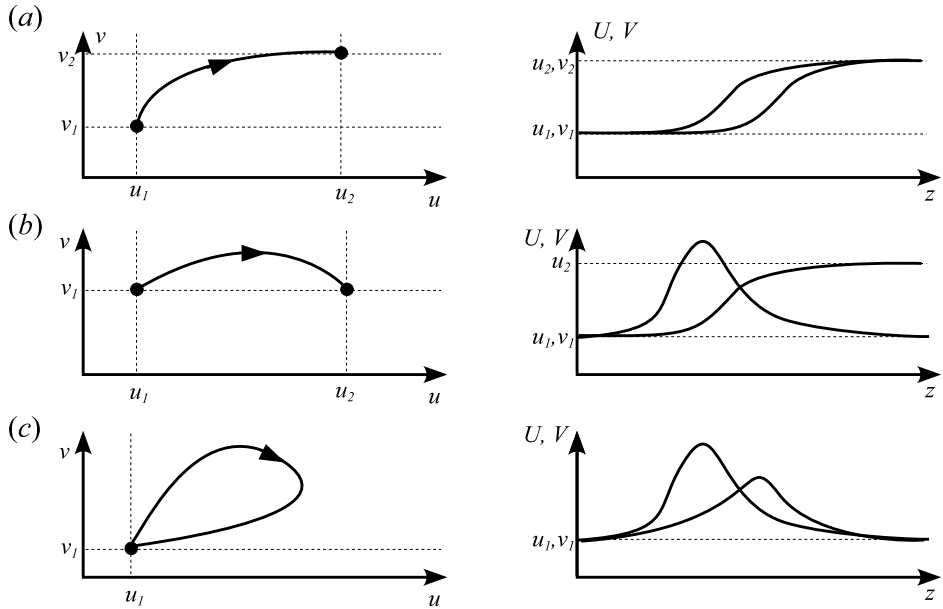}\\
  \caption{Correspondence between bounded traveling wave solutions of system \eqref{basic} (on the right) and the phase curves of the wave system
  \eqref{ws} (on the left); the black dots are singular points of \eqref{ws}. $(a)$ A front-front solution; $(b)$ a front-impulse solution; $(c)$ an impulse-impulse solution}\label{fig:1}
\end{figure}

\begin{proposition}\label{pr1}$ $

\emph{i.} A wave front in $U$ \emph{(}or $V$\emph{)} component
corresponds to a heteroclinic orbit that connects singular points
of \eqref{ws} with different $u$ \emph{(}or $v$\emph{)}
coordinates \emph{(}Fig.\ref{fig:1}a\emph{)};

\emph{ii.} A wave impulse in $U$ \emph{(}or $V$\emph{)} component
corresponds to a heteroclinic orbit that connects singular points
with identical $u$ \emph{(}or $v$\emph{)} coordinates
\emph{(}Fig.\ref{fig:1}b\emph{)} or to a homoclinic curve of a
singular point $(u,\,v)$ of \eqref{ws}
\emph{(}Fig.\ref{fig:1}c\emph{)}.
\end{proposition}

Hereinafter we shall adopt the following terminology: we define
the type of a traveling wave solution of \eqref{basic} with a two
word definition; e.g., a front-impulse solution means that
$u$-profile is a front, and $v$-profile is an impulse (the order
of the terms is important).

For system \eqref{basic} several results on the existence of
one-dimensional traveling waves are known; see, e.g.,
\cite{Horstmann,Keller2,Berez2,Feltman,Berez4,Berez,
   Gueron,   Levine,
Nagai, Stevens}. In most of these references the analysis is
conducted using a particular model which is given in an explicit
form. Quite a different approach was used in \cite{Horstmann}
where the authors consider more general model than \eqref{basic}
and do not restrict themselves to analyzing a model with specific
functions $f(U,V),\,g(U,V)$; instead their aim was to understand
how these functions have to be related to each other in order to
result in traveling wave patterns for $U$ and $V$. We consider a
general class of models as well, and our task is to infer possible
kinds of wave solutions under given restrictions on $f(U,V)$ and
$g(U,V)$.

Our main goal is as follows: we impose some constrains on the
functions $f(U,V),\,g(U,V)$ and study possible traveling wave
solutions with increasing complexity of $f,\,g$. Special attention
is paid to the families of traveling wave solutions such that the
corresponding wave system possesses an infinite number of bounded
orbits. We present the simplest possible models in the form
\eqref{basic} that display traveling wave solutions of a specific
kind.

The main class of the models we deal with is defined in the
following way.

\begin{definition}
We shall call model \eqref{basic} the separable model if
$$
(C1)\qquad f(u,v)=f_1(u)f_2(v),\quad g(u,v)=g_1(u)g_2(v),
$$
where $f_1(u),\,g_1(u)$ are smooth functions for $a\leqslant
u<\infty$; $g_2(v)$ is smooth; $f_2(v)$ is a rational function:
$$
f_2(v)=\frac{Z(v)}{R(v)}\,,
$$
for $b\leqslant v<\infty$; here $a,b>-\infty$ are real constant.

The separable model will be called the reduced separable model if
$(C1)$ holds and
$$
(C2)\qquad f_2(v)g_2(v)\equiv \mbox{const}.
$$
\end{definition}

We organize the paper as follows. In Section \ref{S2} we present
full classification of traveling wave solutions of the reduced
separable models; we also specify necessary and sufficient
conditions for these models to possess specific kinds of traveling
waves. Section \ref{S3} is devoted to the analysis of the
separable model; we show which types of traveling waves can be
expected in addition to the types described in Section \ref{S2};
we also analyze a generalized Keller-Segel model, which does not
belong to the class of the separable models but display a number
of similar properties together with essentially new ones. Section
\ref{S4} contains discussion and conclusions; finally, the details
of numerical computations are presented in Appendix.

\section{Wave solutions of the reduced separable model}\label{S2}
In this section we present the full classification of possible
traveling wave solutions of system \eqref{basic} that satisfies
$(C1),(C2)$. The reason we start with the reduced separable models
is twofold. First, there are models in the literature that have
this particular form (see, e.g., \cite{Levine,  Stevens,Otmer});
second, the special form of the wave system allows the exhaustive
investigation of traveling wave solutions of \eqref{basic}.

The wave system of the reduced separable model reads
\begin{equation}\label{rm}
    \begin{split}
    u' & =cu+f_1(u)g_1(u)/c+\alpha\equiv h(u),\\
    v' & =g_1(u)g_2(v)/c,
\end{split}
\end{equation}
where the first equation is independent of $v$.

\subsection{The cell structure of the phase plane of system \eqref{rm}}

We start with the case of general position. We assume that the
following conditions of non-degeneracy are fulfilled (later we
will relax some of these assumptions):
\begin{equation*}
    \begin{split}
    (A1) & \quad h(u) \mbox{ has no multiple roots};\\
    (A2) & \quad g_2(v) \mbox{ has no multiple roots};\\
    (A3) & \quad g_1(u) \mbox{ and } h(u) \mbox{ have no common roots.}
\end{split}
\end{equation*}

Traveling wave solutions of \eqref{basic} correspond to bounded
orbits of \eqref{rm} different from singular points. Due to the
structure of system \eqref{rm} it is impossible to have a
homoclinic orbit or a limit cycle in the phase plane of
\eqref{rm}, which yields that it is necessary to have at least two
singular points of \eqref{rm} and a heteroclinic orbit connecting
them (see Fig.~\ref{fig:1}a,b) to prove existence of traveling
wave solutions of \eqref{basic} satisfying $(C1),(C2)$.

In general, smooth functions $h(u),\,g(u,v)$ can be written in the
form
\begin{equation*}
    \begin{split}
h(u)&=\tilde{h}(u)(u-\hat{u}_1)\ldots(u-\hat{u}_m),\quad
\tilde{h}(u)\neq0 \mbox{ for any } u,\\
g(u,v)&=g_1(u)\tilde{g}_2(v)(v-\hat{v}_1)\ldots(v-\hat{v}_n),\quad
\tilde{g}_2(v)\neq0 \mbox{ for any } v.
\end{split}
\end{equation*}

The following proposition holds for neighboring roots of $h(u)$
and $g_2(v)$.

\begin{proposition}\label{pr2}$ $

\emph{i.} Let the wave system \eqref{rm} satisfying $(A1)$-$(A3)$
have singular points $(\hat{u},\hat{v}_1)$ and
$(\hat{u},\hat{v}_2)$, where $\hat{v}_1$ and $\hat{v}_2$ are
neighboring roots of $g_2(v)$ then one of these points is a saddle
and the other one is a node.

\emph{ii.} Let the wave system \eqref{rm} satisfying $(A1)$-$(A3)$
have singular points $(\hat{u}_1,\hat{v})$ and
$(\hat{u}_2,\hat{v})$, where $\hat{u}_1$ and $\hat{u}_2$ are
neighboring roots of $h(u)$, then these points can both be
saddles, nodes or one is a node and another is a saddle.
\end{proposition}
\begin{proof}Let $(\hat{u},\hat{v})$ be a singular point of \eqref{rm}. The
eigenvalues of this point $\lambda_1(\hat{u},\hat{v})=h'(\hat{u})$
and $\lambda_2(\hat{u},\hat{v})=g_1(\hat{u})g_2'(\hat{v})$ are
real numbers (henceforth we use prime to denote differentiation
when it is clear with respect to which variable it is carried
out). This implies that singular point $(\hat{u},\hat{v})$ of
system \eqref{rm} cannot be a focus or center.

i. The claim is a simple conjecture of condition $(A2)$.

ii. Let us consider two equilibrium points $(\hat{u}_1,\hat{v})$
and $(\hat{u}_2,\hat{v})$. The eigenvalues
$\lambda_1(\hat{u}_1,\hat{v})$ and $\lambda_1(\hat{u}_2,\hat{v})$
have opposite signs due to $(A1)$.

Consider another pair of eigenvalues
$\lambda_2(\hat{u}_1,\hat{v})$ and $\lambda_2(\hat{u}_2,\hat{v})$
and assume that $(A3)$ holds. If the number of roots of $g_1(u)$
located between $\hat{u}_1$ and $\hat{u}_2$ is even (or zero) then
the signs of these eigenvalues are the same. This implies that one
of the equilibrium points is a saddle whereas the other one is a
node. If the number of roots of $g_1(u)$ located between
$\hat{u}_1$ and $\hat{u}_2$ is odd than the signs of these
eigenvalues are opposite which implies that both equilibriums are
saddles or nodes (one node is attracting and another is
repelling).
\end{proof}

Note, that in case ii of Proposition \ref{pr2} in order to
guarantee that both singular points are nodes one should have
$h'(\hat{u}_1)g_1(\hat{u}_1)g_2'(\hat{v})>0$ and
$h'(\hat{u}_2)g_1(\hat{u}_2)g_2'(\hat{v})>0$. Due to continuity
arguments there exists a family of orbits of \eqref{rm} which tend
to one of the nodes when $z\to\infty$ and to the other node when
$z\to-\infty$.

Taking into account that straight lines $u=\hat{u}_i,\,i=1..m$ and
$v=\hat{v}_j,\,j=1..n$ consist of orbits of system \eqref{rm} we
obtain that the phase plane of \eqref{rm} is divided into
$(m-1)(n-1)$ bounded rectangular domains whose boundaries are
$u=\hat{u}_i,\,i=1..m$ and $v=\hat{v}_j,\,j=1..n$. We shall call
these domains \textit{the orbit cells}.

\begin{figure}[t]
  \centering
  \includegraphics[width=0.8\textwidth]{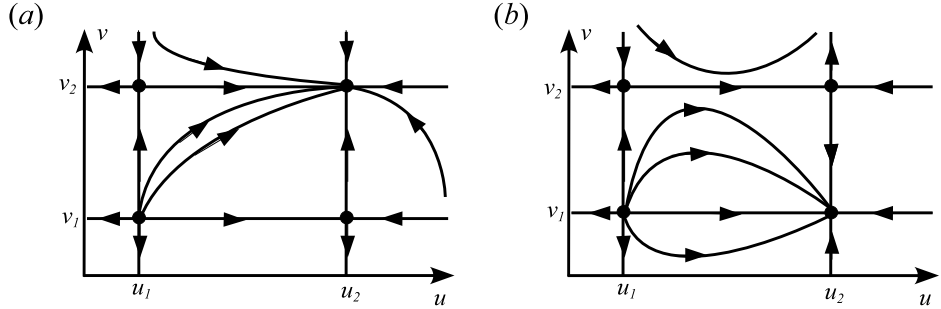}\\
  \caption{Two types of orbit cells of system \eqref{rm}}\label{fig:2}
\end{figure}

Due to Proposition \ref{pr2} it immediately follows that an orbit
cell can be one of the following two types (up to $\pi$-degree
rotation) that are presented in Fig. \ref{fig:2}. The behavior of
the orbits inside a cell is completely described by the types of
the singular points at the corners of the cell. Moreover, any
orbit inside a cell corresponds to a bounded traveling wave
solution of system \eqref{basic}.

Summarizing the previous analysis we obtain the following theorem.

\begin{theorem}\label{th1}
The system \eqref{basic} satisfying $(C1),(C2)$ and $(A1)$-$(A3)$
possesses traveling wave solutions

\emph{i.} of a front-front type \emph{(}Fig. \ref{fig:1}a\emph{)}
if and only if the wave system \eqref{rm} has four singular points
$(u_i,v_i),\,i=1,2$, which are the vertexes of a bounded orbit
cell and every two neighboring vertexes are a node and a saddle
\emph{(}Fig. \ref{fig:2}a\emph{)};

\emph{ii.} of a front-impulse type \emph{(}Fig.
\ref{fig:1}b\emph{)} if and only if the wave system \eqref{rm} has
two neighboring nodes $(\hat{u}_1,\hat{v})$ and
$(\hat{u}_2,\hat{v})$, \emph{(}see Fig. \ref{fig:2}b\emph{)}.
\end{theorem}

\textbf{Remarks to Theorem \ref{th1}.}

1. In both cases the orbits of system \eqref{rm} that correspond
to traveling wave solutions of \eqref{basic} are dense in the
corresponding orbit cell.

2. System \eqref{basic} has a traveling wave solution which is a
front in $v$-component and space-homogeneous in $u$-component if
and only if the wave system \eqref{rm} has neighboring saddle and
node with identical $u$-coordinate, see singular points
$(u_1,v_1), (u_1,v_2)$ and $(u_2,v_1), (u_2,v_2)$ in Fig.
\ref{fig:2}b.

It is possible to write down asymptotics for $u$ and $v$ profiles
(these asymptotics can be used, e.g., as initial conditions for
numerical solutions of \eqref{basic}). We present these
asymptotics only in the simplest case.

Let us assume that the wave system has the form
\begin{equation}\label{as1}
    \begin{split}
    u' & =A(u-\hat{u}_1)(u-\hat{u}_2),\\
      v' &= -B(v-\hat{v}_1)(v-\hat{v}_2)(u-\tilde{u}),
\end{split}
\end{equation}
where $A,B>0$ are constant. An explicit solution of \eqref{as1} is
\begin{equation}\label{as2}
    \begin{split}
    u(z) & =\hat{u}_2+\frac{\hat{u}_1-\hat{u}_2}{1+C_1\exp{\{A(\hat{u}_1-\hat{u}_2)z\}}},\\
     v(z)   &
     =\hat{v}_2+\frac{(\hat{v}_1-\hat{v}_2)(1+C_1\exp{\{A(\hat{u}_1-\hat{u}_2)z\}})^{\frac{B(\hat{v}_2-\hat{v}_1)}{A}}}{(1+C_1\exp{\{A(\hat{u}_1-\hat{u}_2)z\}})^{\frac{B(\hat{v}_2-\hat{v}_1)}{A}}+C_2\exp{\{B(\tilde{u}-\hat{u}_1)(\hat{v}_1-\hat{v}_2)z\}}},
\end{split}
\end{equation}
where $C_1,\,C_2$ are arbitrary constants. We emphasize here that
even with fixed $c$ and $\alpha$ there is a two-parameter family
of wave profiles.

Let $\hat{u}_1<\tilde{u}<\hat{u}_2$ and $\hat{v}_1<\hat{v}_2$. We
consider non-trivial profiles ($C_1C_2\neq 0$). It is
straightforward to show that
$u(-\infty)=\hat{u}_2,u(\infty)=\hat{u}_1$ and, hence, $u$-profile
is a front; $v(-\infty)=\hat{v}_1,v(\infty)=\hat{v}_1$, and
$v$-profile is an impulse. If $\tilde{u}<\hat{u}_1$ or
$\tilde{u}>\hat{u}_2$ then $u$-profile remains the same and
$v$-profile becomes a front.

Formulas \eqref{as2} can be used as a first approximation for wave
profiles of system \eqref{basic} even in the case
$A=A(u),B=B(u,v)$ where $A(u)$ and $B(u,v)$ do not change the sign
when $u\in(\hat{u}_1,\hat{u}_2),\,v\in(\hat{v}_1,\,\hat{v}_2)$.

It is worth noting that if $\hat{u}_1=\tilde{u}$ then
$$
v(\infty)=\hat{v}_2+\frac{\hat{v}_1-\hat{v}_2}{1+C_2},
$$
i.e., the boundary of the profile depends on an arbitrary
constant; we deal with such solutions in the next section.

\subsection{Bifurcations of the travelling wave solutions}
Considering $c$ and $\alpha$ as bifurcation parameters we can
relax some of non-degeneracy conditions $(A1)$-$(A3)$.

First we note that the right-hand side of the second equation of
system \eqref{rm} does not depend in a non-trivial way on $c$ and
$\alpha$ and we will not consider the case when $(A2)$ is
violated. In general, by varying the bifurcation parameters we can
only achieve that either $(A1)$ or $(A3)$ do not hold. We shall
show that in the latter case new traveling wave solutions can
appear in system \eqref{basic} satisfying $(C1),(C2)$.

First let us assume that $(A1)$ does not hold, i.e., function
$h(u)$ has a root $\hat{u}$ of multiplicity $m>1$ for some
$\alpha^*, c^*$, and the wave system \eqref{rm} has a complicated
singular point $(\hat{u},\hat{v})$. The system can be written in
the form
\begin{equation}\label{mulr}
    \begin{split}
    u' & =(u-\hat{u})^m q_1(u),\\
      v'  & =(v-\hat{v})g_1(u)q_2(v),
\end{split}
\end{equation}
where $q_1(\hat{u})q_2(\hat{v})g_1(\hat{u})\neq 0$. Then the
singular point $(\hat{u},\hat{v})$ of system \eqref{mulr} is
either a saddle, a node, or a saddle-node \cite{An}. For the first
two types of critical points, the structure of the phase plane of
the wave system was completely described above. In case of a
saddle-node the line $u=\hat{u}$ divides the plane such that in
one half-plane the singular point is topologically equivalent to a
node, and in the other one it is topologically equivalent to a
saddle; due to the fact that $u=\hat{u}$ consists of solutions of
\eqref{rm} this type of singular points does not yield
qualitatively new bounded solutions of \eqref{rm}. Therefore,
violation of $(A1)$ does not result in new types of wave solutions
of \eqref{basic} satisfying $(C1),(C2)$.

Appearance or disappearance of $u$-fronts correspond to appearance
or disappearance of the roots of the function $h(u)$ which can
occur with variation of the parameters $c$ and $\alpha$. The
simplest case of the appearance of two or three roots corresponds
to the fold or cusp bifurcations respectively \cite{Kuznet} in the
first equation of \eqref{rm}. The simple conditions for the fold
and cusp bifurcations show that, under variation of the boundary
conditions (parameter $\alpha$) and the wave speed (parameter
$c$), appearance of traveling wave solutions of \eqref{basic} is
possible.

\begin{figure}
\centering
\includegraphics[width=0.4\textwidth]{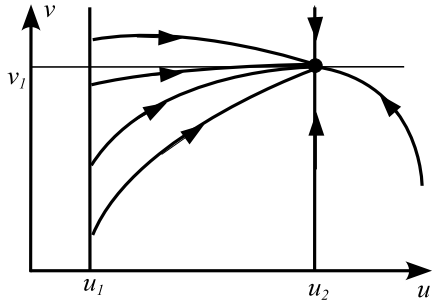}
\caption{The phase plane of system \eqref{rm}; $u_1$ is a root of
both $h(u)$ and $g_1(u)$. The wave solutions corresponding to
bounded orbits of the wave system form a free-boundary
family}\label{fig:FB}
\end{figure}

Now we assume that $(A3)$ does not hold, i.e., the functions
$g_1(u)$ and $h(u)$ have a coinciding root $\hat{u}$. In this case
system \eqref{rm} has a line of non-isolated singularities in the
phase plane $(u,v)$. Each point of the form $(\hat{u},v)$ is a
non-isolated singular point; all the points on the line
$u=\hat{u}$ are either simultaneously attracting or repelling in a
transversal direction to this line \cite{An}. If we assume that
there exists a node $(\hat{u}_1,\hat{v})$ of \eqref{rm} such that
$\hat{u}_1$ is a root of $h(u)$, $g_1(\hat{u}_1)\neq 0$, and there
are no other roots of $h(u)$ between $\hat{u}$ and $\hat{u}_1$
then, due to continuity arguments, there exists a family of
bounded orbits of \eqref{rm} (Fig. \ref{fig:FB}). To describe the
traveling wave solutions corresponding to this family we define
\begin{definition}
We shall say that model \eqref{basic} possesses a family of
free-boundary wave fronts in $v$-component if a\emph{)} every
$v(z)\to\hat{v}$ when $z\to\infty$ $(z\to - \infty)$; b\emph{)}
there exists an interval $(v_1,v_2)$ such that for any
$v^*\in(v_1,v_2)$ it is possible to find $v$-profile with the
property $v(z)\to v^*$ when $z\to-\infty$ $(z\to\infty)$.
\end{definition}

Summarizing we obtain
\begin{proposition}\label{pr3}
The system \eqref{basic} satisfying $(C1),(C2),(A1),(A2)$ has a
traveling wave solution such that $u$-profile is a front and
$v$-profile is a free boundary front if and only if condition
$(A3)$ is violated and there is a node of system \eqref{rm} such
that there are no other singular points of \eqref{rm} between this
node and the line of non-isolated singular points.
\end{proposition}

The primary importance of such traveling wave solutions comes from
the fact that for an arbitrary boundary condition (from a
particular interval) for system \eqref{basic} we can find a wave
solution whose $v$-profile is a front. Note that violation of
$(A3)$ and simultaneous appearance of a free-boundary $v$
component naturally occurs when the roots of $h(u)$ are shifted
under variation of $c$ and $\alpha$. It follows from \eqref{rm}
that bifurcation of $v$-profile occurs at $(c^*,\alpha^*)$ such
that $\hat{u}=-\alpha^*/c^*$ is a simple root of $g_1(u)$.

\subsection{An illustrative example}
Here we present a simple example to illustrate the theoretical
analysis from the previous sections.

We consider model \eqref{basic} with
\begin{equation}\label{fg_ex}
f(u,v)=\frac{(u-l)(1-u)}{v(v-1)},\quad g(u,v)=-k(u-r)v(v-1),
\end{equation}
where $l,k,r$ are non-negative parameters. The particular form of
the functions $f(u,v),g(u,v)$ obviously satisfies $(C1)$ and
$(C2)$.

The wave system reads
\begin{equation}\label{ws_ex}
    \begin{split}
    u' & =cu-k(u-l)(1-u)(u-r)+\alpha,\\
      v'  &= -k(u-r)v(v-1)/c.
\end{split}
\end{equation}
The system \eqref{ws_ex} can have up to six singular points. For
instance, if we fix the parameter values
$l=0.2,r=0.4,k=1,\alpha=-0.1,c=0.3$ then system \eqref{ws_ex}
possesses six singular points; therefore, there are two orbit
cells ensuring existence of traveling wave solutions of system
\eqref{basic}.

\begin{figure}[h]
  \centering
  \includegraphics[width=0.6\textwidth]{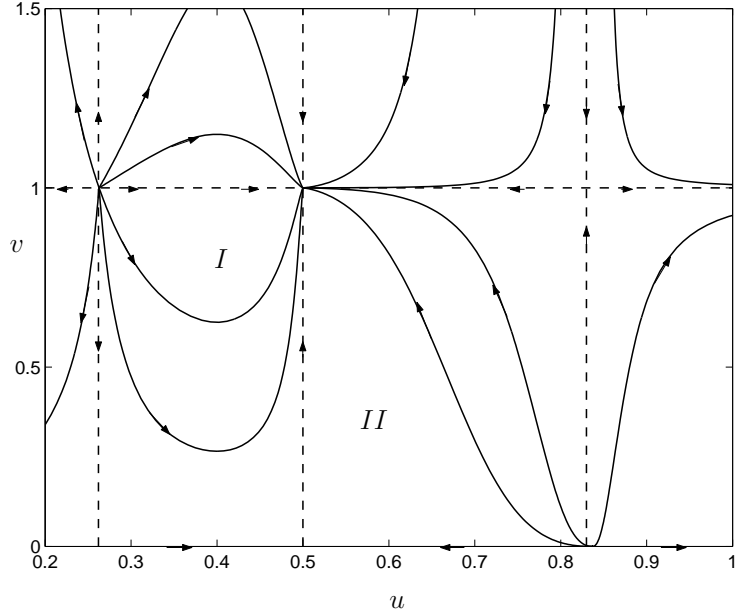}
  \caption{The phase portrait of system \eqref{ws_ex}. The parameters are $l=0.2,\,r=0.4,\,k=1,\,\alpha=-0.1,\,c=0.3$.}\label{fig:3}
\end{figure}

A phase portrait of \eqref{ws_ex} is shown in Fig. \ref{fig:3}.
From Fig. \ref{fig:3} it can be seen that, with the given
parameter values, there exist two qualitatively different
traveling wave solutions of the initial cross-diffusion system
which correspond to two cases of Theorem \ref{th1}.

Numerical solutions of system \eqref{basic} with functions
\eqref{fg_ex} and the given parameter values are shown in Fig.
\ref{fig:4} (the details of the numerical computations are
presented in the Appendix).

If we change the value of $\alpha$ to $-0.12$ then we obtain a
family of free-boundary traveling wave solutions.

\begin{figure}[h!]
  \centering
  \includegraphics[width=0.95\textwidth]{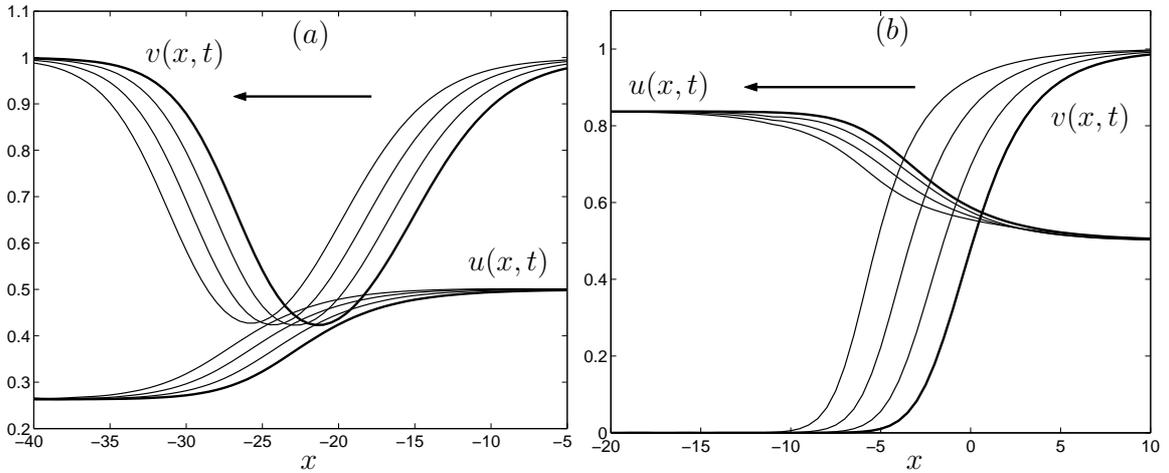}
  \caption{Numerical solutions of system \eqref{basic} with functions \eqref{fg_ex}. The initial conditions are chosen to start the calculations from the orbit cells labelled  $I$ (panel $(a)$, front-impulse solution) or $II$ (panel $(b)$, front-front solution) in Fig. \ref{fig:3}. The solutions are shown for the time moments $t_0=0\mbox{ (bold curves)}<t_1<t_2<t_3=18$ in equal time intervals}\label{fig:4}
\end{figure}

\section{Wave solutions of the separable models and some generalizations} \label{S3}
\subsection{General theory}
In this section we study models \eqref{basic} which satisfy
$(C1)$. The rational function $f_2(v)=Z(v)/R(v)$ can be presented
in the form
$$
f_2(v)=\frac{Z(v)}{R(v)}=\frac{\tilde{Z}(v)(v-v_1)\ldots(v-v_m)}{\tilde{R}(v)(v-\breve{v}_1)\ldots(v-\breve{v}_k)},
$$
where $\tilde{Z}(v),\tilde{R}(v)$ do not have real roots;
$m\geqslant 0,k>0$; $v_i\neq\breve{v}_j$ for any $i,j$. The wave
system has the form
\begin{equation}\label{rmod}
    \begin{split}
    u' & =cu+\alpha+f_1(u)g_1(u)g_2(v)\frac{Z(v)}{cR(v)},\\
     v'   &=g_1(u)g_2(v)/c.
\end{split}
\end{equation}
By transformation of the independent variable
\begin{equation}\label{f:11}
    dy=\frac{dz}{cR(v)},
\end{equation}
which is smooth for any $v$ except for $v=\breve{v}_i$, $i=1..k$,
system \eqref{rmod} becomes
\begin{equation}\label{rmod1}
    \begin{split}
    \frac{du}{dy} & =c^2R(v)(u+\alpha/c)+f_1(u)g_1(u)g_2(v)Z(v),\\
     \frac{dv}{dy}   &=g_1(u)g_2(v)R(v).
\end{split}
\end{equation}
The roots of functions $f_1(u),g_1(u),g_2(v),Z(v),R(v)$ do not
depend on parameters $c$ and $\alpha$, hence, we will suppose that
the following conditions of non-degeneracy are fulfilled:
\begin{equation*}
    \begin{split}
    (B1) & \quad R(v) \mbox{ and } g_2(v) \mbox{ have no common
    roots;}\\
    (B2) & \quad f_1(u) \mbox{ and } g_1(u) \mbox{ have no common roots;}\\
    (B3) & \quad f_1(u), g_1(u), g_2(v), R(v) \mbox{ have no multiple roots.}
\end{split}
\end{equation*}
Coordinates of singular points $(u,v)=(\hat{u},\hat{v})$ of
\eqref{rmod1} can be found from one of the systems:
\begin{align}
    \label{sys1} -\alpha/c&=u,\quad g_2(v)=0,\\
    \label{sys2} R(v)&=0,\quad f_1(u)=0,\\
    \label{sys3} R(v)&=0,\quad g_1(u)=0,
\end{align}
or from combination of \eqref{sys1}-\eqref{sys3}.

To infer possible types of the singular points of \eqref{rmod1} we
consider $D=\Tr^2 J-4\Det J$, where $J$ is the Jacobian of
\eqref{rmod1} evaluated at a singular point $(\hat{u},\hat{v})$.

If $(\hat{u},\hat{v})$ is a solution of \eqref{sys1} then
\begin{align*}
\Tr J&=R(\hat{v})(c^2+g_1(\hat{u})g_2'(\hat{v})),\\
\Det J&=c^2g_1(\hat{u})R(\hat{v})^2g_2'(\hat{v}),\\
D&=R(\hat{v})^2(c^2-g_1(\hat{u})g_2'(\hat{v}))^2;
\end{align*}
if $(\hat{u},\hat{v})$ is a solution of \eqref{sys2} then
\begin{align*}
\Tr
J&=g_1(\hat{u})g_2(\hat{v})(Z(\hat{v})f_1'(\hat{u})+R'(\hat{v})),\\
\Det
J&=(g_1(\hat{u})g_2(\hat{v}))^2Z(\hat{v})f_1'(\hat{u})R'(\hat{v}),\\
D&=(g_1(\hat{u})g_2(\hat{v}))^2(Z(\hat{v})f_1'(\hat{u})-R'(\hat{v}))^2;
\end{align*}
if $(\hat{u},\hat{v})$ is a solution of \eqref{sys3} then
$$
\Det J=0,\quad \Tr
J=f_1(\hat{u})g_2(\hat{v})Z(\hat{v})g_1'(\hat{u}).
$$

Consequently we obtain that $(\hat{u},\hat{v})$ is a saddle or a
node for the cases corresponding to \eqref{sys1} and \eqref{sys2},
and is a saddle, node, or saddle-node (see \cite{An}) in the case
\eqref{sys3}. Just as for the reduced separable models there are
no singular points of \eqref{rmod1} of center or focus type.

Here we do not pursue the problem of classification of possible
structures of the phase plane of \eqref{rmod1} and are only
concerned with new types of traveling wave solutions. Analyzing
formulas  for $\Det J$ and $\Tr J$ and applying arguments in the
line with the proof of Proposition \ref{pr2} and Theorem \ref{th1}
we obtain

\begin{proposition}\label{pr4}
Let coordinates of singular points of \eqref{rmod1} satisfy
\eqref{sys1} and function $g_2(v)$ have two real neighboring roots
$\hat{v}_1$ and $\hat{v}_2$. If the function $R(v)$ has an odd
number of roots between $\hat{v}_1$ and $\hat{v}_2$ and point
$(\hat{u}=-\alpha/c,\hat{v}_1)$ is a node, then
$(\hat{u}=-\alpha/c,\hat{v}_2)$ is a node as well \emph{(}and vice
versa\emph{)}.

Let coordinates of singular points of \eqref{rmod1} satisfy
\eqref{sys2} and $u=\hat{u}$ be a root of $f_1(u)$,  $R(v)$ have
two real neighboring roots $\hat{v}_1$ and $\hat{v}_2$. If the
function $Z(v)$ has an odd number of roots between $\hat{v}_1$ and
$\hat{v}_2$ and point $(\hat{u},\hat{v}_1)$ is a node, then
$(\hat{u},\hat{v}_2)$ is a node as well \emph{(}and vice
versa\emph{)}.
\end{proposition}

Due to the structure of system \eqref{rmod1} the phase plane is
divided into horizontal strips, whose boundaries are given by
$v=\hat{v}$, where $\hat{v}$ is a root of $R(v)$ or $g_2(v)$; all
singular points of \eqref{rmod1} are situated on these boundaries.
Bringing in the continuity arguments we obtain that under the
conditions of Proposition \ref{pr4} there is a family of bounded
orbits of \eqref{rmod1} which correspond to the traveling wave
solution of \eqref{basic} of an impulse-front type.

\begin{figure}
\centering
\includegraphics[width=0.95\textwidth]{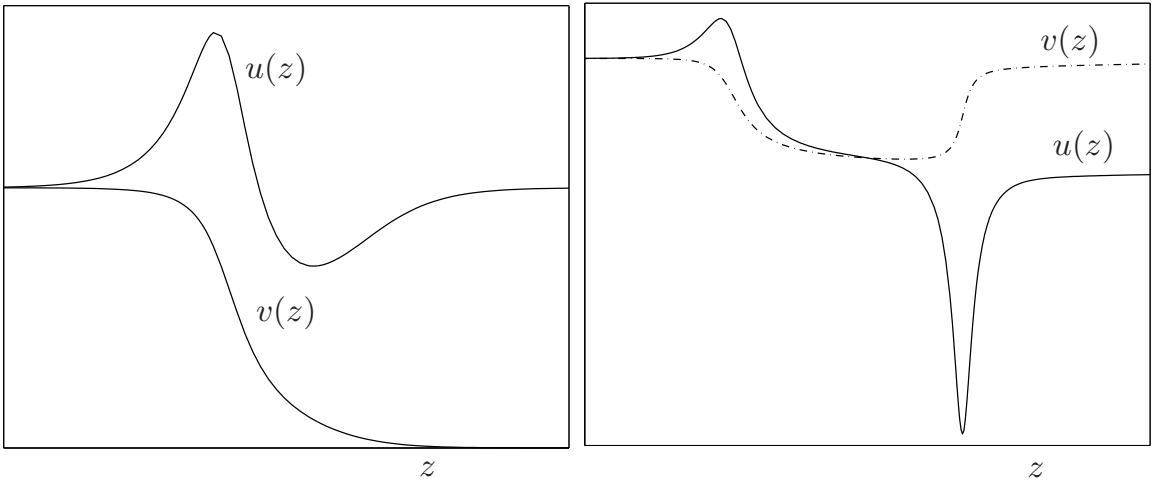}
\caption{Complex shape of wave profiles of model \eqref{rmod}. The
functions are
$f_1(u)=u-1,\,g_1(u)=u,\,g_2(v)=1,\,R(v)=v(v-1),\,Z(v)=2v-1$. The
left panel shows an impulse-front solution, the right panel shows
a front-impulse solution}\label{fig:NN}
\end{figure}

It is worth noting that the structure of the phase plane inside a
strip can be quite arbitrary, and we can only indicate asymptotic
behavior of orbits in neighborhoods of singular points. As a
result, under given boundary conditions (or, equivalently, fixed
$\alpha$) families of traveling wave solutions may have complex
shapes (opposite to the examples presented in Fig. \ref{fig:1}).
For instance, there can be non-monotonous fronts with humps or
impulses which also have multiple humps and hollows. In general,
we can only state that the form of impulses and fronts can be
quite arbitrary, which is illustrated in Fig. \ref{fig:NN}.

Under variation of parameters $c$ and $\alpha$ it is possible that
function $g_1(u)$ has a root $\hat{u}=-\alpha^*/c^*$ for some
values of the parameters; in this case system \eqref{rmod1} has a
line of non-isolated singular points $u=\hat{u}$ and the analysis
in this situation is similar to the analysis which led to
Proposition \ref{pr3}: a family of free-boundary fronts appears in
$v$ component.

Now let us relax the condition $(B2)$; we assume that there exists
such $\hat{u}$ that $f_1(\hat{u})=0$ and $g_1(\hat{u})=0$. We can always find values of
$c$ and $\alpha$ such that $\hat{u}=-\alpha/c$. In this case the
phase plane of \eqref{rmod1} has a line $u=\hat{u}$ of
non-isolated singular points. After the change of the independent
variable $d\tau=dy/(u-\hat{u})$ the resulting system still
possesses singular point of the form $(\hat{u},\hat{v})$, where
$\hat{v}$ is a root of $R(v)$. If this point is a node then,
applying continuity arguments, we obtain that there exist a family
of orbits of system \eqref{rmod1} such that some solutions from a
neighborhood of $(\hat{u},\hat{v})$ tend to this point if
$y\to\infty$ (or $y\to-\infty$) and tend to point of the form
$(\hat{u},v^*)$ if $y\to-\infty$ (or $y\to\infty$), where $v^*$ is
an arbitrary constant from some interval. These solutions
correspond to traveling wave solution of \eqref{basic} such that
$u$-profile is an impulse and $v$ component is a free-boundary
front (see Fig. \ref{fig:6}).

Summarizing we obtain
\begin{theorem}\label{th2}
The system \eqref{basic} satisfying $(C1)$ and $(B1)$-$(B3)$ can
only possess traveling wave solutions of the following kinds:

\emph{i.} front-front solutions;

\emph{ii.} front-impulse solutions;

\emph{iii.} impulse-front solutions;

\emph{iv.} Under variation of parameters $c$ and $\alpha$ it is
possible to have wave solutions where $u$ component is a front,
$v$ component is a free-boundary front;

\emph{v.} Under the additional condition that $(B2)$ does not hold
it is possible to have wave solutions with $u$ component is an
impulse and $v$ component is a free-boundary front.
\end{theorem}

\subsection{The phase plane analysis of the Keller-Segel model}
The classical Keller-Segel model has the form \eqref{basic} with
$f(u,v)=\delta u/v$, $g(u,v)=-ku$, where $\delta,k>0$ (see
\cite{Keller2}). In our terminology this model falls in the class
of the separable models; the wave system reads
\begin{equation}\label{ks1}
    \begin{split}
    u' &=cu-\delta k \frac{u^2}{cv}+\alpha,\\
     v'   & =-ku/c,
\end{split}
\end{equation}
which, with the help of transformation \eqref{f:11}, can be
reduced to the form \eqref{rmod1}:
\begin{equation}\label{ks2}
    \begin{split}
    \frac{du}{dy} &=c^2uv-\delta k u^2+\alpha c v,\\
     \frac{dv}{dy}   & =-kuv.
\end{split}
\end{equation}

\begin{figure}[b!]
\centering
\includegraphics[width=0.7\textwidth]{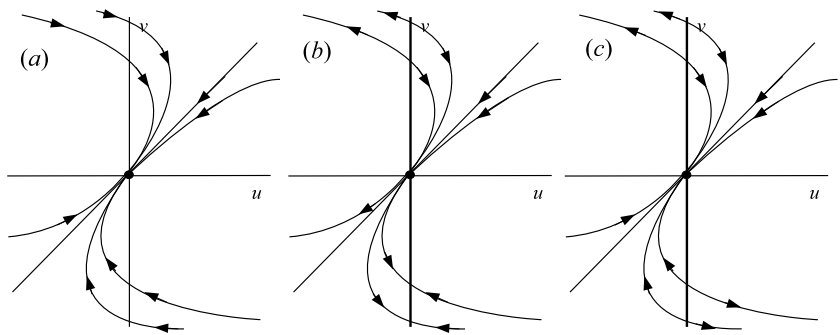}
\caption{The phase planes of systems \eqref{ks3} $(a)$,
\eqref{ks2} $(b)$ and \eqref{ks1} $(c)$}\label{fig:6s}
\end{figure}

If $\alpha\neq 0$ then system \eqref{ks2} has the only singular
point $(u,v)=(0,0)$. In this case the separable model cannot
possess traveling wave solutions (Theorem \ref{th2}). Hence we put
$\alpha=0$. Note that the last requirement is necessary if one
supposes that $U(x,t)$ should be finite.

For $\alpha=0$ system \eqref{ks2} has a line of non-isolated
singular points $(0,v)$, and there is also additional degeneracy
at the point $(0,0)$ (case v. in Theorem \eqref{th2}). This can be
seen applying the second transformation of the independent
variable $d\tau=udy$, which leads to the system
\begin{equation}\label{ks3}
    \begin{split}
    \frac{du}{d\tau} &=c^2v-\delta k u,\\
     \frac{du}{d\tau}   & =-kv,
\end{split}
\end{equation}
for which the origin is a topological node with the eigenvalues
$\lambda_1=-k\delta$ and $\lambda_2=-k$.  Thus there exists a
family of traveling wave solutions whose $u$-profile is an impulse
and $v$-profile is a free-boundary front.

In Fig. \ref{fig:6s} we show how the parametrization of the phase
curves of the wave system change after the transformations of the
independent variables. This picture can also serve as an
illustration to assertion v. of Theorem \ref{th2}.

Due to biological interpretation of the Keller-Segel model it is
necessary to have $U(x,t)\geqslant 0$ and $V(x,t)\geqslant 0$. Using the fact that
an eigenvector $(1,(\delta-1)k/c^2)$ corresponds to $\lambda_1$
and $(1,0)$ corresponds to $\lambda_2$ it is straightforward to
see that to ensure the existence of non-negative travelling way
solutions we should have $\delta>1$.

\begin{figure}[h!]
\includegraphics[width=0.95\textwidth]{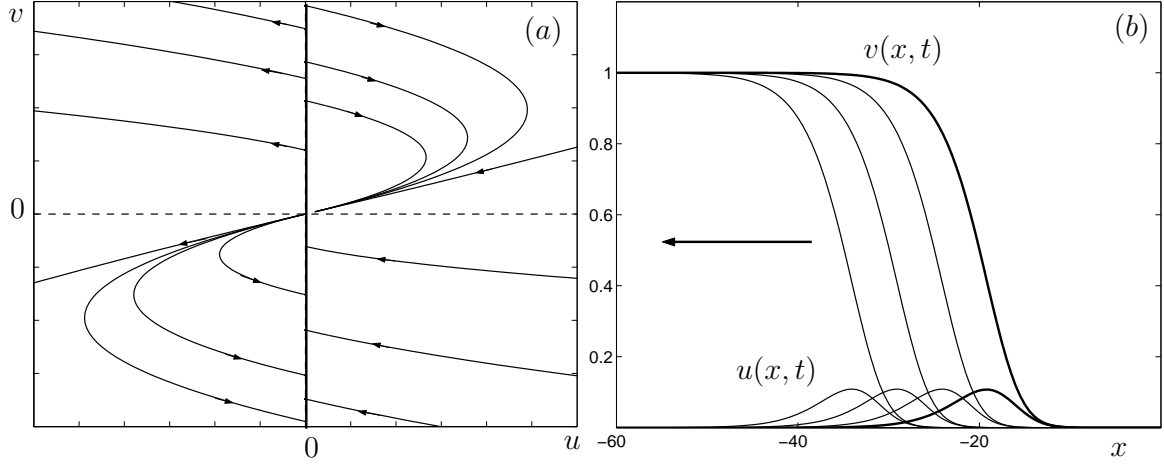}
\caption{$(a)$ The phase plane of system \eqref{ks1}; the
parameters are $\alpha=0, k=1,\delta=4,c=0.5$. $(b)$ Numerical
solutions of the Keller-Segel system with the parameters given in
$(a)$. The solutions are shown for the time moments $t_0=0\mbox{
(bold curves)}<t_1<t_2<t_3=30$ in equal time
intervals}\label{fig:6}
\end{figure}

Numerical solutions of the Keller-Segel model are given in Fig.
\ref{fig:6}b. Originally, the Keller-Segel model was suggested to
describe movement of bands of \textit{E. coli} which were observed
to travel at a constant speed when the bacteria are placed in one
end of a capillary tube containing oxygen and an energy source
\cite{Keller2}. In Fig. \ref{fig:6}b it can be seen that bacteria
($u(x,t)$) seek an optimal environment: the bacteria avoid low
concentrations and move preferentially toward higher
concentrations of some critical substrate ($v(x,t)$). The
stability of the traveling solutions found was studied
analytically in \cite{Gueron,{Nagai}}.

\subsection{Impulse-impulse solutions}
In the preceding sections we studied the systems \eqref{basic} where the functions $f(u,v),g(u,v)$ can be
represented as a product of functions that only depend on one
variable. The next natural step is to assume that these functions
depend on affine expressions $au+bv+c$, where $a,b,c$ are not
equal to zero simultaneously. Here we present an explicit example
of such a system. The example is motivated by appearance of a
particular type of traveling wave solutions, which is absent in
the separable models.

We suppose that
\begin{equation}\label{rd1}
    \begin{split}
    f(u,v) &=\frac{\delta u}{\beta u+v},\quad \delta>0,\beta\geqslant 0,\\
     g(u,v)&=-ku+rv,\quad k,r\geqslant 0.
\end{split}
\end{equation}
If $\beta,r=0$ then we have the Keller-Segel model studied in
Section 3.2.

The wave system for \eqref{basic} with the functions given by
\eqref{rd1} reads
\begin{equation}\label{rd2}
    \begin{split}
    u' &=cu+\frac{\delta u(-ku+rv)}{c(\beta u+v)},\\
      v'  &=(-ku+rv)/c,
\end{split}
\end{equation}
where we put parameter $\alpha$ equal zero.

After the change of the independent variable $dz/(c(\beta
u+v))=dy$ system \eqref{rd2} takes the form
\begin{equation}\label{rd3}
    \begin{split}
    \frac{du}{dy} &=c^2u(\beta u+v)+\delta u(-ku+rv),\\
     \frac{dv}{dy}  &=(-ku+rv)(\beta u+v).
\end{split}
\end{equation}
If $r=0$ then system \eqref{rd3} has a line on non-isolated
singular points $u=0$; if $r\neq 0$ then $(u,v)=(0,0)$ is an
isolated non-hyperbolic singular point of \eqref{rd3} (i.e., both
eigenvalues of the Jacobian evaluated at this point are zero).

First we consider the case $r=0$. After yet another transformation
$d\tau=dy/u$ we obtain
\begin{equation}\label{rd4}
    \begin{split}
    \frac{du}{d\tau} &=(c^2\beta-\delta k)u+c^2v,\\
      \frac{dv}{d\tau}  &=-k(\beta u+v).
\end{split}
\end{equation}
Thus for $\delta>1$ and $0\leqslant\beta\leqslant
k(\sqrt{\delta}-1)^2/c^2$ the origin is a node for system
\eqref{rd4} which implies that in the initial system \eqref{rd3}
there exist a family of bounded orbits which tent to $(0,0)$ for
$y\to\infty$ or $y\to-\infty$. This family corresponds to a family
of traveling wave solution of the system \eqref{basic} where
$u$-profile is an impulse and $v$-profile is a free-boundary
front. The picture is topologically equivalent to the phase
portrait shown in Fig. \ref{fig:6}a.
\begin{figure}[t!]
\includegraphics[width=0.95\textwidth]{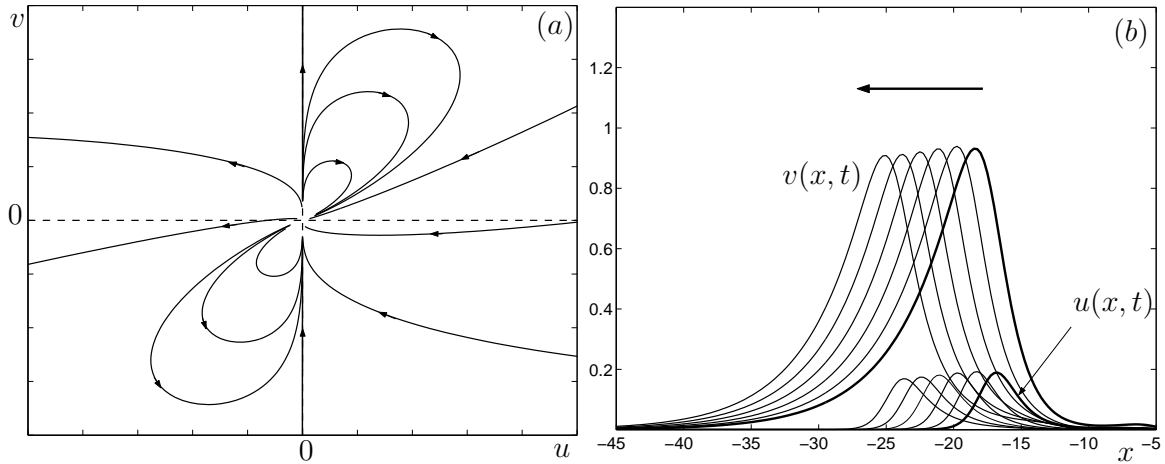}
\caption{$(a)$ The phase plane of system \eqref{rd2}; the
parameters are $k=1,\beta=1,\delta=2,c=0.43,r=0.1$. $(b)$
Numerical solutions of system \eqref{basic} with the functions
given by \eqref{rd1} and the parameter values as in $(a)$. The
solutions are shown for the time moments $t_0=0\mbox{ (bold
curves)}<t_1<t_2<t_3<t_4<t_5=20$ in equal time
intervals}\label{fig:7}
\end{figure}

For $r\neq 0$ the wave system \eqref{rd3} has singular point
$(0,0)$ possessing two elliptic sectors in its neighborhood (see
Fig. \ref{fig:7}a). The proof of existence of the elliptic sectors
can easily be conducted with the methods given in \cite{Berez5}.
Asymptotics of homoclinics composing the elliptic sector are $u=0$
(trivial) and $v=K^+u$, where $K^+$ is the biggest root of the
equation
$$
K^2(c^2+r(\delta-1))+K(\beta c^2-k(\delta-1)-\beta r)+\beta k=0.
$$
The family of homoclinics in the phase plane $(u,v)$ correspond to
the family of wave impulses for the system \eqref{basic} (see Fig.
\ref{fig:1}c and Fig. \ref{fig:7}b). To our knowledge such kind of
solutions (infinitely many traveling wave solutions of
impulse-impulse type with the fixed values of the parameters) was
not previously described in the literature.

The results of the numerical computations (Fig. \ref{fig:7}b)
indicate that the family of traveling impulses is clearly
non-stationary, since its amplitude decays visibly in time. Which
is important, however, is that it is possible to observe moving
impulses at least at a finite time interval.

\section{Conclusions}\label{S4}

In this paper we described all possible traveling wave solutions
$\{u(z),\, v(z),\, z=x+c\,t\}$ of the cross-diffusion
two-component PDE model \eqref{basic} satisfying $(C1)$, where the
cross-diffusion coefficient may depend on both variables and
possess singularities. Such kind of models is widely used in
modeling populations that can chemotactically react to an
immovable signal (attractant) (e.g., \cite{Horstmann,Keller2} and
references therein). The study of traveling wave solutions of
model \eqref{basic} was carried out by qualitative and bifurcation
analysis of the phase portraits of the wave system \eqref{ws} that
depends on parameters $c$ and $\alpha$. Here $c$ is a speed of
wave propagation and $\alpha$ characterizes boundary conditions
under given $c$.  Any traveling wave solution of \eqref{basic}
with given boundary conditions corresponds to a solution of wave
system \eqref{ws} with specific values of  parameters $c$ and
$\alpha$; the converse is also true. Therefore, instead of trying
to construct a traveling wave solution with the given boundary
conditions, we study the set of all possible bounded solutions of
the wave system considering $c$ and $\alpha$ as its parameters.
This approach allows identifying all boundary conditions for which
model \eqref{basic} possesses traveling wave solutions.

The main attention is paid to the so-called separable model, i.e.,
to model \eqref{basic} that satisfies $(C1)$; it is worth noting
that in this case the functions $f(u,v)$ and $g(u,v)$ in
\eqref{basic} are products of factors that depend on a single
variable. We showed that for some fixed values of parameters $c$
and $\alpha$ the solutions of the wave system compose
two-parameter family (for explicit formulas in a simple particular
case see \eqref{as2}). One result is that $\{u(z),
v(z)\}$-profiles of the wave solution of the reduced separable
models \eqref{basic} that satisfies $(C1),(C2)$ can be only
front-front or front-impulse; for more general case of the
separable models we can additionally have impulse-front solutions.
For some special relations between $c$, $\alpha$ and the model
parameters model \eqref{basic} can have wave solutions whose
$v$-profiles are free-boundary fronts, i.e., $v(z)$ tends to an
arbitrary constant from some interval at $z\to\infty$ or at
$z\to-\infty$. Note that traveling wave solutions of the
well-known Keller-Segel model \cite{Keller2} as well as of its
generalization (Section 3.3) have impulse-front profiles with a
free-boundary front (Sections 3.2, 3.3).

We also considered a natural extension of the separable models;
namely, we gave an example of model \eqref{basic} where $f,g$ are
products of factors that depend on affine expression of both
variables (see \eqref{rd1}). This model can be considered as a
generalization of the Keller-Segel model because it has two
additional parameters and turns into the Keller-Segel model if
both of these parameters are zero. If only one of the parameters
vanishes, the model has a family of ``Keller-Segel''-type
solutions, i.e., $u$-profile is an impulse and $v$-profile is a
frond with a free boundary. Importantly, in some parameter domains
this model possesses a two-parameter family of impulse-impulse
solutions (Fig. \ref{fig:7}). To the best of our knowledge, such
type of traveling wave solutions was not previously described in
the literature: depending on the initial conditions traveling
impulses can have quite a different form for the fixed values of
the system parameters. Taking into account the fact that such
solutions are absent in the separable models, we can consider
model \eqref{basic} with functions \eqref{rd1} as the simplest
model possessing this type of traveling wave solutions.

Rearrangements of traveling wave solutions of PDE model
\eqref{basic}, which occur with changes of the wave propagation
velocity and the boundary conditions, correspond to bifurcations
of its ODE wave system. In particular, appearance/disappearance of
front-profiles with variation of parameters $c$ and $\alpha$
correspond to the fold or cusp bifurcations in the wave system;
rearrangement of a front to an impulse can be accompanied by
appearance/disappearance of non-isolated singular points in the
phase plane of the corresponding wave system (see Section 2.3).
Existence of non-isolated singular points in the wave system may
result in the existence of free-boundary fronts in model
\eqref{basic}. For instance, this is the case for the Keller-Segel
model.

We emphasize that the separable model \eqref{basic}, when the
values of parameters $c$ and $\alpha$ fixed in the wave system,
possesses, in general, two-parameter family of traveling wave
solutions. There are infinitely many bounded orbits of \eqref{ws}
that correspond to traveling wave solutions of \eqref{basic} (see
Figures \ref{fig:2}, \ref{fig:4}, \ref{fig:6}, \ref{fig:7}). It is
of particular interest that in all numerical solutions of
\eqref{basic} that we conducted it is possible to observe
traveling waves. We did not discuss the issue of stability of the
traveling wave solutions found, but we note that it is usually
true that unstable solutions cannot be produced in numerical
calculations. It is tempting to put forward a hypothesis that the
presence of additional degrees of freedom (two free parameters) is
the reason of producing traveling waves in numerical computations.
This important question can be a subject of future research.

\section{Appendix}\label{Appendix}

We did numerical simulations of system \eqref{basic} for
$x\in[-L,\,L]$, where $L$ varied in different numerical
experiments. We used no-flux boundary conditions for the spatial
variable. Inasmuch as we wanted to study the behavior of the
traveling wave solutions in an infinite space we chose such space
interval so that to avoid the influence of boundaries.

We used an explicit difference scheme. The approximation of the
taxis term is an "upwind" explicit scheme \cite{Morton} which is
frequently used for cross-diffusion systems (e.g.,
\cite{Tsyganov}). More precisely,
\begin{equation*}
\begin{split}
    u_i^{t+1} & =u_i^t+\frac{\Delta t}{(\Delta x)^2}(u_{i+1}^t-2u_i^t+u_{i-1}^t)-\frac{\Delta t}{(\Delta x)^2}(A(v_{i+1}^t-v_i^t)-B(v_{i}^t-v_{i-1}^t)),\\
     v_i^{t+1}&= v_i^t+(\Delta t) g(u_i^t,v_i^t), \quad
     i=2,...,N-1,
\end{split}
\end{equation*}
where for the positive taxis (pursuit) (i.e., $f(u,v)<0$),
\begin{equation*}
\begin{split}
    A & =f(u_i^t,v_i^t)\quad \mbox{if } v_{i+1}\geqslant v_i, \\
    A & =f(u_{i+1}^t,v_{i+1}^t)\quad \mbox{if } v_{i+1}<v_i,\\
    B & =f(u_{i-1}^t,v_{i-1}^t)\quad \mbox{if } v_{i}\geqslant v_{i-1}, \\
    B & =f(u_{i}^t,v_{i}^t)\quad \mbox{if } v_{i}<v_{i-1}.
\end{split}
\end{equation*}
For the negative taxis (invasion):
\begin{equation*}
\begin{split}
    A & =f(u_i^t,v_i^t)\quad \mbox{if } v_{i+1}< v_i, \\
    A & =f(u_{i+1}^t,v_{i+1}^t)\quad \mbox{if } v_{i+1}\geqslant v_i,\\
    B & =f(u_{i-1}^t,v_{i-1}^t)\quad \mbox{if } v_{i}<v_{i-1}, \\
    B & =f(u_{i}^t,v_{i}^t)\quad \mbox{if } v_{i}\geqslant v_{i-1}.
\end{split}
\end{equation*}
We used $\Delta t=0.001$, $\Delta x=0.1$. For the boundary
conditions:
\begin{equation*}
\begin{split}
    u_1^{t+1}&=u_2^t, \quad u_N^{t+1}=u_{N-1}^t,\\
    v_1^{t+1}&=v_2^t, \quad v_N^{t+1}=v_{N-1}^t.
\end{split}
\end{equation*}
For the initial conditions we used numerical solutions of the
corresponding wave systems.

\paragraph{Acknowledgments.} FSB and GPK express their gratitude to Dr. A. Stevens for numerous useful discussions on the first version of the present work. The work of FSB has been supported in part by NSF Grant
\#634156 to Howard University.


\begin{thebibliography}{100}

\bibitem{Horstmann} Horstmann, D., and A. Stevens. 2004. A constructive approach to traveling waves in chemotaxis. Journal of Nonlinear Science 14(1):1-25.



\bibitem{Keller1} E.F. Keller, L.A. Segel,  Model for chemotaxis, Journal of Theoretical Biology, 30(2) (1971) 225-234.

\bibitem{Keller2} E.F. Keller, L.A. Segel, Traveling Bands of Chemotactic Bacteria - Theoretical Analysis, Journal of Theoretical Biology, 30(2) (1971)
235-248.

\bibitem{Erban} R. Erban, H.G. Othmer, Taxis equations for amoeboid cells, J Math Biol, 54(6) (2007) 847-885.

\bibitem{Ni} W.-M. Ni, Diffusion, cross-diffusion and their spike-layer steady states, Not. Am. Math. Soc. 45(1) (1998)
9-18.

\bibitem{Bud} E.O. Budrene, H.C. Berg, Complex Patterns Formed by Motile Cells of Escherichia-Coli, Nature 349(6310) (1991)
630-633.

\bibitem{Mikhailov} A.A. Samarskii, A.P. Mikhailov, Principles of mathematical modeling: ideas, methods, examples. Taylor \&
Francis, 2002.

\bibitem{Berez2} F.S. Berezovskaya, G.P. Karev, Bifurcations of travelling waves in population taxis models, Uspekhi Fizicheskikh Nauk 169(9) (1999) 1011-1024.

\bibitem{Feltman} D.L. Feltham, M.A.J. Chaplain, Travelling waves in a model of species migration, Applied Mathematics Letters 13(7)(2000) 67-73.

\bibitem{Habib} S. Habib, C. Molina-Paris, T.S. Deisboeck, Complex dynamics of tumors: modeling an emerging brain tumor system with coupled reaction-diffusion equations, Physica A: Statistical Mechanics and its Applications 327(3-4) (2003) 501-524.

\bibitem{Sherratt} J.A Sherratt. Traveling wave solutions of a mathematical model for tumor encapsulation, Siam Journal on Applied Mathematics 60(2) (1999) 392-407.

\bibitem{tsyganov} M.A. Tsyganov, V.N.Biktashev, J. Brindley, A.V. Holden, G.R. Ivanitsky, Waves in systems with cross-diffusion as a new class of nonlinear waves, Uspekhi Fizicheskikh Nauk 177(3) (2007)
275-300.

\bibitem{Berez4} F.S. Berezovskaya, G.P. Karev, Traveling waves in cross-diffusion models of the dynamics of populations, Biofizika 45(4) (2000) 751-756.

\bibitem{An} A.A. Andronov, E.A. Leontovich, I.I. Gordon, A.G. Maier. Qualitative Theory of Second-order Dynamic Systems, John Wiley \&
Sons, 1973.

\bibitem{Kuznet} Y.A. Kuznetsov, Elements of Applied Bifurcation Theory
Springer, 2004.

\bibitem{Berez} F. Berezovskaya, G. Karev, Parametric portraits of travelling waves of population models with polynomial growth and auto-taxis rates, Nonlinear Analysis: Real World Applications 1(1) (2000) 123-136.

\bibitem{Murray} J.D. Murray, Mathematical Biology, Springer,
2005.

\bibitem{Volpert} A.I. Volpert, V.A. Volpert, V.A. Volpert, Traveling Wave Solutions of Parabolic Systems,
AMS, 1994.

\bibitem{Gueron} S. Gueron, N. Liron, A Model of Herd Grazing as a Traveling Wave, Chemotaxis and Stability, Journal of Mathematical Biology 27(5) (1989) 595-608.

\bibitem{Levine} H.A. Levine, B.D. Sleeman, A system of reaction diffusion equations arising in the theory of reinforced random walks, Siam Journal on Applied Mathematics 57(3) (1997) 683-730.

\bibitem{Nagai} T. Nagai, T. Ikeda, Traveling Waves in a Chemotactic Model, Journal of Mathematical Biology 30(2) (1991) 169-184.

\bibitem{Stevens} A. Stevens, Trail following and aggregation of myxobacteria, J. Biol. Sys. 3(4) (1995) 1059-1068.

\bibitem{Otmer} H.G. Othmer, A. Stevens, Aggregation, blowup, and collapse: The ABC's of taxis in reinforced random walks, Siam Journal on Applied Mathematics 57(4) (1997) 1044-1081.

\bibitem{Berez5} F.S. Berezovskaya, A.S. Novozhilov, G.P. Karev, Population models with singular equilibrium, Mathematical Biosciences 208 (2007) 270-299.

\bibitem{Morton} K. W. Morton, D. F. Mayers, Numerical solutions
of partial differential equations, Cambrige University Press,
1989.

\bibitem{Tsyganov}
M. A. Tsyganov, J. Brindley, A. V. Holden, V. N. Biktashev,
Soliton-like phenomena in one-dimensional cross-diffusion systems:
a predator-prey pursuit and evasion example, Physica D: Nonlinear
Phenomena, 197(1-2) (2004) 18-33.

\end{thebibliography}
\end{document}